\documentclass{amsart}
\usepackage{cite}
\usepackage{amsmath,amssymb,amsfonts}
\usepackage{algorithmic}
\usepackage{graphicx}
\usepackage{textcomp}

\usepackage{epsfig}
\usepackage{latexsym}
\usepackage{pdfsync}

\newcommand{\eq}{\begin{equation}\begin{array}{rllllllllllllllllllllllllllllllll}}
\newcommand{\ee}{\end{array}\end{equation}}
\newcommand{\bmt}{\left[ \begin{array}{ccccccccc}}
\newcommand{\emt}{\end{array}\right]}
\newcommand{\bea}{\begin{eqnarray}}
\newcommand{\eea}{\end{eqnarray}}
\newcommand{\bean}{\begin{eqnarray*}}
\newcommand{\eean}{\end{eqnarray*}}
\newcommand{\bc}{\begin{center}}
\newcommand{\ec}{\end{center}}

\title{Boundary Control of the Beam Equation by Linear Quadratic Regulation }
\thanks{This work was supported by AFOSR under FA9550-20-1-0318 }

\author{Arthur J. Krener}
\address{Naval Postgraduate School\\ 833 Dyer Road Monterey\\ CA 93943-5216}
\email{ajkrener@nps.edu}


\begin{document}

\begin{abstract}
We present and solve a Linear Quadratic Regulator (LQR) for the boundary control of the beam equation.  We use the simple technique of completing the square to get an explicit solution.
By decoupling the spatial frequencies we are able to reduce an infinite dimensional LQR to an infinte family of two two dimensional LQRs each of which can be solved explicitly.
\end{abstract}

\maketitle

\section{Introduction}
We consider the stabilization of the linear beam equation using  a Linear Quadratic Regulator (LQR).   By decoupling the spatial frquencies  we obtain a complete
and explicit solution to the LQR including the closed loop eigenvalues.  The simple technique of completing the square yields  Riccati PDEs.   At each spatial frequency
the Riccati PDEs reduce to the algebraic Riccati equation of a two dimensional problem which is readily solvable.
The sums of the optimal cost and optimal feedback of these two dimensional problems yield the optimal cost and optimal feedback of the infinite dimensional LQR.
The only technical issues that arise are whether these sums are convergent.  We discuss which Lagrangian yield convergence.

The study of optimal control of systems governed by partial differential equations goes back at least to Lions \cite{JL71}.  More recent treatises on this topic 
are the works of Curtain and Zwart \cite{CZ95}, \cite{CZ20}, Lasiecka and Trigiani \cite{LT00} and Krstic and Smyshlyaev \cite{KS08} who use backstepping.    
LQR boundary control has been used by Lasiecka and Trigiani \cite{LT86}, Burns and King \cite{BK95}, Hulsing \cite{Hu99}, Burns and Hulsing \cite{BH01},
Cristofaro, DeLuca and Lanari \cite{CDL}.  Coron, D'Andrea and Bastin \cite{CAB07} found Lyapunov functions for the boundary control of hyperbolic conservation laws. Guo et al. have considered boundary control of the beam equation in the presence of disturbaces,
\cite{GZA14}, \cite{GK14}.  Other papers on boundary control of the beam equation are Morgul \cite{Mo02}, Militec and Arnold \cite{MA13}, Han, Li, Xu and Liu \cite{HXL16}.

More recently we introduced the Completing the Square technique to solve LQR problems for partial differential equations.  We solved an  LQR problem
for the heat equation under distributed control in \cite{Kr20a} and  under boundary control in \cite{Kr20b}.  In both cases using an extension of Al'brekht's 
method \cite{Al61} we were able to find the higher degree terms in the Taylor polynomial expansions of the optimal cost and the optimal feedback.

In \cite{Kr21} we solved an LQR for the boundary control of the wave equation by decoupling the spatial frequencies.  This allowed us to reduce an infinite 
dimensional LQR problem to an infinite family of two dimensional LQR problems each of which can be explicitly solved.  In this paper we show that an LQR problem 
for the beam equation can also be reduced to an infinite family of two dimensional LQR problems each of which can be explicitly solved.

\section{Boundary Control of the Beam Equation}
In Exercise 3.18 Curtain and Zwart consider the undamped   beam equation
subject to boundary control action
\bean
\frac{\partial^2 f}{\partial t^2}(x,t)&=& -\frac{\partial^4 f}{\partial x^4}(x,t)\\
f(0,t)=0,&& f(1,t)=0\\
\frac{\partial^2 f}{\partial x^2}(0,t)=0,&&\frac{\partial^2 f}{\partial x^2}(1,t)=u(t)
\\
f(x,0)=f_1(x),&& \frac{\partial f}{\partial t}(x,0)=f_2(x)
\eean

It is convenient to intoduce vector notation, let
\bean
z(x,t)&=& \bmt z_1(x,t)\\  z_2(x,t) \emt \ = \ \bmt f(x,t)\\ \frac{\partial f}{\partial t}(x,t)\emt
\eean
We also allow damping 
\bean
{d\over dt} z(x,t)&=&{\it A}z(x,t) 
\eean
where $\alpha\ge0$ and ${\it A}$ is the matrix differential operator
\bean
{\it A}&=& \bmt 0&1\\ -\frac{\partial^4 }{\partial x^4}& -\alpha \emt,
\eean
The boundary conditions are
\bean
z_1(0,t)=0,&& z_1(1,t)=0\\
\frac{\partial^2 z_1}{\partial x^2}(0,t)=0,&&\frac{\partial^2 z_1}{\partial x^2}(1,t)=u(t)
\eean
and the initial conditions are 
\bean
z_1(x,0)=f_1(x),&&z_2 (x,0)=f_2(x)
\eean
If $\alpha=0$ the beam is undamped and if $\alpha>0$ the beam is damped.

The eigenvalues of the open loop system are
\bean
\lambda_n&=&{-\alpha+ \mbox{ sign}(n)\sqrt{\alpha^2-4n^4\pi^4}     \over 2}
\eean
and the corresponding eigenvectors are
\bean
v_n&=& \bmt {1\over \lambda_n}\\ 1\emt \sin |n| \pi x
\eean
for $n=\pm 1,\pm 2,\pm 3,\ldots$. 
If $|n|$ is small enough the eigenvalues $\lambda_n, \lambda_{-n}$ can be real numbers but
for large $|n|$ the eigenvalues $\lambda_n, \lambda_{-n}$ are complex and conjugate.

We wish to find a feedback to stabilize the beam,  $z(x,t)\to 0$ as $t\to \infty$.
Another possibility is that we wish to stabilize the beam  to some   open loop trajectory $z^*(x,t)$.
We define $\tilde{z}(x,t)=z(x,t)-z^*(x,t)$ and we seek a feedback to drive   
$\tilde{z}(x,t)\to 0$ as $t\to \infty$. Because of linearity these are  equivalent problems so we only consider the first one.

 We shall use
a Linear Quadratic Regulator (LQR) to find the desired feedback.  We choose a $2\times 2 $ nonnegative definite
matrix valued function $Q(x_1,x_2)$ which is symmetric in $x_1,x_2$, $Q(x_1,x_2)=Q(x_2,x_1)$,
and a positive scalar $R$.  Consider the problem of minimizing 
\bea \label{crit}
\int_0^\infty \left(\iint_{\it S}  z'(x_1,t)Q(x_1,x_2) z(x_2,t)\  dA + Ru^2(t)\right) \ dt
\eea
subject to the beam dynamics where ${\it S}=[0,1]^2$ and $dA=dx_1dx_2$.

Let $P(x_1,x_2)$ be any $2\times 2 $ symmetric 
matrix valued function  which is also symmetric in $x_1,x_2$, $P(x_1,x_2)=P(x_2,x_1)$.
Suppose there is a control trajectory $u(t)$ such that  the corresponding state trajectory $z(x,t)$ goes to $0$ as $t \to \infty$ then by the Fundamental
Theorem of Calculus 
\bean
&&0=\iint_{\it S}  z'(x_1,t)P(x_1,x_2) z(x_2,t)\  dA \\
&&+\int_0^\infty \iint_{\it S} {d \over dt} \left(z'(x_1,t)P(x_1,x_2) z(x_2,t)\right)\  dA  \ dt
\eean
We expand the integrand of the time integral.
\bean
&&0=\iint_{\it S}  z'(x_1,t)P(x_1,x_2) z(x_2,t)\  dA \\
&&+\int_0^\infty \iint_{\it S} {d z'\over dt}(x_1,t) P(x_1,x_2) z(x_2,t)\  dA  \ dt\\
&&+\int_0^\infty \iint_{\it S}  z'(x_1,t)Q(x_1,x_2) {d z\over dt}(x_2,t)\  dA  \ dt
\eean
\bean
&&0=\iint_{\it S}  z'(x_1,t)P(x_1,x_2) z(x_2,t)\  dA \\
&&+\int_0^\infty \iint_{\it S} \bmt z_2(x,t)\\-\frac{\partial^4 z_1}{\partial x^4}(x_,t)-\alpha z_2(x_1,t)\emt'\bmt  P_{1,1}(x_1,x_2)& P_{1,2}(x_1,x_2) \\P_{2,1}(x_1,x_2)&P_{2,2}(x_1,x_2)\emt \bmt z_1(x_2,t)\\z_2(x_2,t)\emt\  dA  \ dt\\
&&+\int_0^\infty \iint_{\it S}  \bmt z_1(x_2,t)\\z_2(x_2,t)\emt'\bmt  P_{1,1}(x_1,x_2)& P_{1,2}(x_1,x_2) \\
P_{2,1}(x_1,x_2)&P_{2,2}(x_1,x_2)\emt  \bmt z_2(x,t)\\-\frac{\partial^4 z_1}{\partial x_2^4}(x,t)-\alpha z_2(x_2,t)\emt \  dA  \ dt
\eean

\bea \label{FTC}
&&0=\iint_{\it S}  z'(x_1,t)P(x_1,x_2) z(x_2,t)\  dA \\
&&+\int_0^\infty \iint_{\it S} z_2(x_1,t) P_{1,1}(x_1,x_2)z_1(x_2,t) +z_1(x_1,t) P_{1,1}(x_1,x_2)z_2(x_2,t)   \nonumber \\
&&+ z_2(x_1,t) P_{1,2}(x_1,x_2)z_2(x_2,t)+z_2(x_1,t) P_{2,1}(x_1,x_2)z_2(x_2,t) \nonumber\\
&&  -\frac{\partial^4 z_1}{\partial x_1^4}(x_1,t)P_{2,1}(x_1,x_2)z_1(x_2,t)  -z_1(x_1,t)P_{1,2}(x_1,x_2)  \frac{\partial^4 z_1}{\partial x_2^4}(x_2,t) \nonumber \\
&&  -\frac{\partial^4 z_1}{\partial x_1^4}(x_1,t)P_{2,2}(x_1,x_2)z_2(x_2,t)  -z_2(x_1,t)P_{2,2}(x_1,x_2)  \frac{\partial^4 z_1}{\partial x_2^4}(x_2,t) 
\\
&&-\alpha z_2(x_1,t)P_{2,1}(x_1,x_2)z_1(x_2,t)-\alpha z_2(x_1,t)P_{2,2}(x_1,x_2)z_2(x_2,t)\nonumber\\
&&-\alpha z_1(x_1,t)P_{1,2}(x_1,x_2)z_2(x_1,t)-\alpha z_2(x_1,t)P_{2,2}(x_1,x_2)z_2(x_1,t)
 \nonumber \ dA \ dt  \nonumber
\eea

We assume that $ P(x_1,x_2)$ satisfies these boundary conditions
\bea
P(0,x_2)= P(1,x_2)=P(x_1,0)=P(x_1,1)=0 \label{Pbc1}\\
\frac{\partial^2 P}{\partial x_1}(0,x_2)=\frac{\partial^2 P}{\partial x_1}(1,x_2)=\frac{\partial^2 P}{\partial x_2}(x_1,0)=\frac{\partial^2 P}{\partial x_2}(x_1,1)=0
 \label{Pbc2}
\eea

We integrate by parts four times to get these identities
\bean
&&\iint_{\it S }-\frac{\partial^4 z_1}{\partial x_1^4}(x_1,t)P_{2,1}(x_1,x_2)z_1(x_2,t)\ dA\\
&& =\iint_{\it S }u(t)\frac{\partial P_{2,1}}{\partial x_1}(1,x_2)z_1(x_2,t)-z_1(x_1,t)\frac{\partial^4 P_{2,1}}{\partial x_1^4}(x_1,x_2)z_1(x_2,t)\ dA\\
\\
  && \iint_{\it S }-z_1(x_1,t)P_{1,2}(x_1,x_2)  \frac{\partial^4 z_1}{\partial x_2^4}(x_2,t)\  dA\\
  && =\iint_{\it S }-z_1(x_1,t)\frac{\partial^4 P_{1,2}}{\partial x_2^4}(x_1,x_2)z_1(x_2,t)+z_1(x_1,t)\frac{\partial P_{1,2}}{\partial x_2}(x_1,1)u(t)\ dA\\
  \\
  &&\iint_{\it S }-\frac{\partial^4 z_1}{\partial x_1^4}(x_1,t)P_{2,2}(x_1,x_2)z_2(x_2,t)\ dA\\
&& =\iint_{\it S }u(t)\frac{\partial P_{2,2}}{\partial x_1}(1,x_2)z_2(x_2,t)-z_1(x_1,t)\frac{\partial^4 P_{2,2}}{\partial x_1^4}(x_1,x_2)z_2(x_2,t)\ dA\\
\\
  && \iint_{\it S }-z_1(x_1,t)P_{2,2}(x_1,x_2)  \frac{\partial^4 z_1}{\partial x_2^4}(x_2,t)\  dA\\
  && =\iint_{\it S }-z_2(x_1,t)\frac{\partial^4 P_{2,2}}{\partial x_2^4}(x_1,x_2)z_1(x_2,t)+z_1(x_1,t)\frac{\partial P_{1,2}}{\partial x_2}(x_1,1)u(t)\ dA
\eean

We plug these identities into (\ref{FTC} ) and obtain 
\bea \label{FTC1}
&&0=\iint_{\it S}  z'(x_1,t)P(x_1,x_2) z(x_2,t)\  dA \\
&&+\int_0^\infty \iint_{\it S} z_2(x_1,t) P_{1,1}(x_1,x_2)z_1(x_2,t) +z_1(x_1,t) P_{1,1}(x_1,x_2)z_2(x_2,t) \nonumber\\
&&+ z_2(x_1,t) P_{1,2}(x_1,x_2)z_2(x_2,t)+z_2(x_1,t) P_{2,1}(x_1,x_2)z_2(x_2,t) \nonumber\\
&&+u(t)\frac{\partial P_{2,1}}{\partial x_1}(1,x_2)z_1(x_2,x_2)-z_1(x_1,t)\frac{\partial^4 P_{2,1}}{\partial x_1^4}(x_1,x_2)z_1(x_2,t)\  \nonumber\\
&&+z_1(x_1,t)\frac{\partial^4 P_{1,2}}{\partial x_2^4}(x_1,x_2)z_1(x_2,t)+z_1(x_1,t)\frac{\partial P_{1,2}}{\partial x_2}(x_1,1)u(t)\  \nonumber\\
&& +u(t)\frac{\partial P_{2,2}}{\partial x_1}(1,x_2)z_2(x_2,t)-z_1(x_1,t)\frac{\partial^4 P_{2,2}}{\partial x_1^4}(x_1,x_2)z_2(x_2,t)\  \nonumber\\
&& -z_2(x_1,t)\frac{\partial^4 P_{2,2}}{\partial x_2^4}(x_1,x_2)z_1(x_2,t)+z_1(x_1,t)\frac{\partial P_{1,2}}{\partial x_2}(x_1,1)u(t) \nonumber
\\
&&-\alpha z_2(x_1,t)P_{2,1}(x_1,x_2)z_1(x_2,t)-\alpha z_2(x_1,t)P_{2,2}(x_1,x_2)z_2(x_2,t)\nonumber\\
&&-\alpha z_1(x_1,t)P_{1,2}(x_1,x_2)z_2(x_1,t)-\alpha z_2(x_1,t)P_{2,2}(x_1,x_2)z_2(x_1,t)
 \nonumber\ dA\ dt  \nonumber
\eea

We add the the right side of (\ref{FTC1}) to the criterion (\ref{crit}) to get an equivalent criterion

\bea \label{crit1}
&&\int_0^\infty \left(\iint_{\it S}  z'(x_1,t)Q(x_1,x_2) z(x_2,t)\  dA + Ru^2(t)\right) \ dt\\
&&\iint_{\it S}  z'(x_1,t)P(x_1,x_2) z(x_2,t)\  dA \nonumber \\
&&+\int_0^\infty \iint_{\it S} z_2(x_1,t) P_{1,1}(x_1,x_2)z_1(x_2,t) +z_1(x_1,t) P_{1,1}(x_1,x_2)z_2(x_2,t) \nonumber\\
&&+ z_2(x_1,t) P_{1,2}(x_1,x_2)z_2(x_2,t)+z_2(x_1,t) P_{2,1}(x_1,x_2)z_2(x_2,t) \nonumber\\
&&+u(t)\frac{\partial P_{2,1}}{\partial x_1}(1,x_2)z_1(x_2,x_2)-z_1(x_1,t)\frac{\partial^4 P_{2,1}}{\partial x_1^4}(x_1,x_2)z_1(x_2,t)\  \nonumber\\
&&+z_1(x_1,t)\frac{\partial^4 P_{1,2}}{\partial x_2^4}(x_1,x_2)z_1(x_2,t)+z_1(x_1,t)\frac{\partial P_{1,2}}{\partial x_2}(x_1,1)u(t)\  \nonumber\\
&& +u(t)\frac{\partial P_{2,2}}{\partial x_1}(1,x_2)z_2(x_2,t)-z_1(x_1,t)\frac{\partial^4 P_{2,2}}{\partial x_1^4}(x_1,x_2)z_2(x_2,t)\  \nonumber\\
&& -z_2(x_1,t)\frac{\partial^4 P_{2,2}}{\partial x_2^4}(x_1,x_2)z_1(x_2,t)+z_1(x_1,t)\frac{\partial P_{1,2}}{\partial x_2}(x_1,1)u(t) \nonumber
\\
&&-\alpha z_2(x_1,t)P_{2,1}(x_1,x_2)z_1(x_2,t)-\alpha z_2(x_1,t)P_{2,2}(x_1,x_2)z_2(x_2,t)\nonumber\\
&&-\alpha z_1(x_1,t)P_{1,2}(x_1,x_2)z_2(x_1,t)-\alpha z_2(x_1,t)P_{2,2}(x_1,x_2)z_2(x_1,t)
 \nonumber\ dA\ dt  \nonumber
\eea

We would like to find a $2\times 1$ matrix valued function 
\bea \label{K}
K(x)&=& \bmt K_1(x) & K_2(x)\emt
\eea
so that the time integrand in (\ref{crit1}) is a perfect square of the form
\bean
\iint_{\it S}  \left(u(t) -K(x_1)z(x_1,t)\right)R \left(u(t) -K(x_2)z(x_2,t)\right) \ dA
\eean
Clearly the terms quadratic in $u(t)$ agree so we equate terms containing 
the product of $u(t)$ and $z(x_2,t)$.   This yields the equation
\bean
-R \bmt K_1(x_2)& K_2(x_2)\emt&=& \bmt \frac{\partial P_{2,1}}{\partial x_1}(1,x_2)& \frac{\partial P_{2,2}}{\partial x_1}(1,x_2)\emt
\eean
so we set
\bea \label{K1}
\bmt K_1(x_2) & K_2(x_2)\emt&=&-R^{-1} \bmt \frac{\partial P_{2,1}}{\partial x_1}(1,x_2)& \frac{\partial P_{2,2}}{\partial x_1}(1,x_2)\emt
\eea
By symmetry
\bean
\bmt K_1(x_1) & K_2(x_1)\emt&=&-R^{-1} \bmt \frac{\partial P_{1,2}}{\partial x_2}(x_1,1)& \frac{\partial P_{2,2}}{\partial x_2}(x_1,1)\emt
\eean

Then we equate terms containing the product of $z(x_1,t)$ and $z(x_2,t)$ and we obtain the Riccati PDEs for the boundary control
of the beam equation,
\bea \label{RicPDEs11}
&&\frac{\partial^4 P_{1,2}}{\partial x_2^4}(x_1,x_2)+\frac{\partial^4 P_{2,1}}{\partial x_1^4}(x_1,x_2)
+Q_{1,1}(x_1,x_2)\\
&&=\gamma^2 \frac{\partial P_{1,2}}{\partial x_2}(x_1,1)\frac{\partial P_{2,1}}{\partial x_1}(1,x_2) \nonumber\\  \nonumber\\
\label{RicPDEs12}
&&P_{1,1}(x_1,x_2)-\frac{\partial^4 P_{2,2}}{\partial x_1^4}(x_1,x_2) 
-\alpha z_1(x_1,t)P_{1,2}(x_1,x_2)z_2(x_1,t)+Q_{1,2}(x_1,x_2) \\
&&=\gamma^2 \frac{\partial P_{1,2}}{\partial x_2}(x_1,1)\frac{\partial P_{2,2}}{\partial x_1}(1,x_2)\nonumber\\ \nonumber\\
\label{RicPDEs21}
&&P_{1,1}(x_1,x_2)-\frac{\partial^4 P_{2,2}}{\partial x_2^4}(x_1,x_2)
-\alpha z_2(x_1,t)P_{2,1}(x_1,x_2)z_1(x_2,t) +Q_{2,1}(x_1,x_2)\\
&&=\gamma^2 \frac{\partial P_{2,2}}{\partial x_2}(x_1,1)\frac{\partial P_{2,1}}{\partial x_1}(1,x_2) \nonumber\\ \nonumber\\
\label{RicPDEs22}
&&P_{1,2}(x_1,x_2)+P_{2,1}(x_1,x_2) -2\alpha z_2(x_1,t)P_{2,2}(x_1,x_2)z_2(x_2,t)+Q_{2,2}(x_1,x_2)\\
&&=\gamma^2 \frac{\partial P_{2,2}}{\partial x_2}(x_1,1)\frac{\partial P_{2,2}}{\partial x_1}(1,x_2) \nonumber
\eea
where $\gamma^2 =R^{-1}\beta^2$.

To simplify the problem we decouple the spatial frequencies by assuming that $Q(x_1,x_2)$ has the expansion
\bea \label{Qsum}
Q(x_1,x_2)&=&  \sum_{n=1}^\infty  \bmt Q^{n,n}_{1,1}&Q^{n,n}_{1,2}\\Q^{n,n}_{2,1}&Q^{n,n}_{2,2}\emt \sin n\pi x_1 \sin n\pi x_2
\eea
and  $Q^{n,n}_{1,2}=Q^{n,n}_{2,1}$.

We assume that $P(x_1,x_2)$ has a similar expansion
\bea  \label{Psum}
P(x_1,x_2)&=& \sum_{n=1}^\infty \bmt P^{n,n}_{1,1}&P^{n,n}_{1,2}\\P^{n,n}_{2,1}&P^{n,n}_{2,2}\emt \sin n\pi x_1  \sin n\pi x_2
\eea
with $P^{m,n}_{1,2}=P^{m,n}_{2,1}$.
Clearly any such $P(x_1,x_2)$ satisfies the boundary conditions (\ref{Pbc1}) and (\ref{Pbc2}).

Then (\ref{K1}) implies 
\bea \label{Ksum}
K(x_2)&=&-R^{-1}\beta  \sum_{n=1}^\infty \bmt P_{2,1}^{m,n} &P_{2,2}^{m,n} \emt n\pi \sin n\pi x 
\eea
and the Riccati PDEs   (\ref{RicPDEs11}, \ref{RicPDEs12}, \ref{RicPDEs21}, \ref{RicPDEs22}) imply that 
\bea
0&=&   -2n^4\pi^4P^{n,n}_{1,2} +Q^{n,n}_{1,1}-n^2\pi^2\gamma^2\left(P^{n,n}_{1,2}\right)^2  \label{11}\\
0&=& P^{n,n}_{1,1} -n^4\pi^4 P^{n,n}_{2,2}-\alpha P^{n,n}_{1,2}  +Q^{n,n}_{1,2} -n^2\pi^2\gamma^2P^{n,n}_{1,2}P^{n,n}_{2,2} \label{12}\\
0&=& P^{n,n}_{1,1} -n^4\pi^4 P^{n,n}_{2,2}-\alpha P^{n,n}_{2,1} +Q^{n,n}_{2,1} -n^2\pi^2\gamma^2P^{n,n}_{2,2}P^{n,n}_{2,1}\label{21}\\
0&=& 2P^{n,n}_{1,2}-2\alpha P^{n,n}_{2,2} +Q^{n,n}_{2,2}-n^2\pi^2 \gamma^2\left(P^{n,n}_{2,2}\right)^2 \label{22}
\eea
where $\gamma^2= R^{-1}\beta^2$.

For each $n=1,2,\ldots$ these are the Riccati equations of  the two dimensional  LQR with matrices 
\eq  \label{lqr2}
F^{n,n}=\bmt 0&1\\ -n^4 \pi^4 &-\alpha \emt, && G^{n,n}=\bmt 0\\n\pi \beta\emt\\
Q^{n,n}=\bmt Q^{n,n}_{1,1}& Q^{n,n}_{1,2}\\Q^{n,n}_{2,1}&Q^{n,n}_{2,2}\emt,&& R^{n,n}=\bmt R\emt
\ee

We use the quadratic formula to solve   (\ref{11})
\bea \label{p12}
P^{n,n}_{1,2}&=& {-n^2\pi^2\pm \sqrt{n^4\pi^4+{\gamma^2Q^{n,n}_{1,1}\over n^2\pi^2}}\over \gamma^2}
\eea
then   (\ref{22})  implies
\bea \label{p22}
P^{n,n}_{2,2}&=&{-\alpha\pm \sqrt{\alpha^2+ n^2\pi^2 \gamma^2\left(Q^{n,n}_{2,2}+2P^{n,n}_{1,2}\right)}\over n^2\pi^2\gamma^2}
\eea
Since we want $P^{n,n}_{2,2}$ to be nonnegative we  take the positive sign.
Then  (\ref{12}) implies 
\bea \label{p11}
P^{n,n}_{1,1}&=& \alpha P^{n,n}_{1,2}+ n^4\pi^4 P^{n,n}_{2,2}-Q^{n,n}_{1,2}+n^2\pi^2\gamma^2P^{n,n}_{1,2}P^{n,n}_{2,2}
\eea
If the two dimensional LQR (\ref{lqr2}) satisfies the standard conditions then the associated Riccati equation has a unique nonnegative definite solution.  This implies that if we take the negative sign in (\ref{p12}) the resulting $P^{n,n}$ is not nonnegative definite. 


The $2\times 2$ closed loop system is 
\bean
F^{n,n}+GK^{n,n}&=& \bmt 0&1\\ -n^4\pi^4 -\gamma^2 P_{2,1}^{n,n}& -\alpha-n^2\pi^2\gamma^2 P_{2,2}^{n,n}\emt
\eean
and the closed loop eigenvalues are
\bea \nonumber
\mu_n&=& 
-{\alpha +n^2\pi^2\gamma^2P^{n,n}_{2,2}\over 2}+\left(\mbox{sign n}\right){\sqrt{(\alpha+n^2\pi^2 \gamma^2P^{n,n}_{2,2})^2-4(n^4\pi^4+\gamma^2 P^{n,n}_{2,1})}\over 2}\\  \label{cleva}
\eea
for $n=\pm 1,\pm2, \pm3,\ldots$.
For $n=1,2,\ldots$ the corresponding eigenvectors of the $n^{th}$ $2\times 2$ closed loop system are 
 \bea \label{eigv2}
 \bmt {1\over \mu_n}\\1 \emt
 ,&\quad\ & \bmt {1\over \mu_{-n}}\\1 \emt
 \eea
 The corresponding eigenvectors of the infinite dimensional closed loop system are 
 \bea \label{eigv2in}
v_n(x)= \bmt {1\over \mu_n}\\ 1\emt \sin |n|\pi x  ,&\quad\ &  v_{-n}(x)= \bmt {1\over \mu_{-n}}\\ 1\emt \sin |n|\pi x
 \eea
 Notice that at least for large $|n|$, $\mu_n$ and $\mu_{-n}$ are complex conjugates as are $v_n(x)$ and $v_{-n}(x)$.
 
The trajectories of the infinite closed loop system are 
\bean
z(x,t)&=& \sum_{n=-\infty}^\infty \zeta_n(t) \bmt {1\over \mu_n}\\1 \emt \sin |n|\pi x
\eean
where
\bean
\zeta_n(t) &=& e^{\mu_n t}  \zeta_n^0
\eean
If  $\mu_n$ and $\mu_{-n}$ are complex and conjugatethen  $\zeta_n^0$ and $\zeta_{-n}^0$ must be complex conjugates for $z(x,t)$ to be
real valued.

Notice we can control each spatial frequency independently.   If we don't want to damp out the $n^{th} $ spatial frequency then we set
$Q^{n,n}=0$ so that  $P^{n,n}=0$ and $K^{n,n}=0$. 

But we must address the questions of whether (\ref{Qsum}) and (\ref{Psum})  converge.
If there is an $N>0$ and an $r>1$ such that 
\bean
\left\| Q_{i,j}^{n,n}\right\|_\infty&\le& {q\over n^r}
\eean
for $i,j=1,2$ and $n>N$ 
then clearly (\ref{Qsum}) converges.  We assume that we have chosen $Q^{n,n} $ such that this is true.

We apply the Mean Value Theorem to  (\ref{p12}) to obtain
\bean 
P^{n,n}_{1,2}&=& {1\over 2  s^{1/2}} {Q^{n,n}_{1,1}\over n^2\pi^2}
\eean
for some $s$ between $n^2\pi^2$ and $\sqrt{n^4\pi^4+{\gamma^2Q^{n,n}_{1,1}\over n^2\pi^2}}$.
Since ${1\over 2  s^{1/2}}$ is monotonically decreasing on this interval and takes on its maximum
value at $n^2\pi^2$ we conclude that
\bean 
P^{n,n}_{1,2}&\le&{Q^{n,n}_{1,1}\over 2n^3\pi^3}\ \le \ {q\over 2n^{3+r}\pi^3}
\eean
so clearly the sum
\bean
P_{1,2}(x_1,x_2)&=& P_{2,1}(x_1,x_2)\ =\ \sum_{n=0}^\infty P^{n,n}_{1,2}\sin nx_1 \sin nx_2
\eean
converges.

If $\alpha=0$ then (\ref{p22}) implies that
\bean 
P^{n,n}_{2,2}&=&{1\over  n\pi \gamma}\sqrt{Q^{n,n}_{2,2}+2P^{n,n}_{1,2}}\ \le\ {c\over n^{1+r/2}}
\eean
for $n>N$ and some constant $c$
so clearly the sum
\bea\label{P22sum}
P_{2,2}(x_1,x_2)&=&  \sum_{n=0}^\infty P^{n,n}_{2,2}\sin nx_1 \sin nx_2
\eea 
converges.

If $\alpha>0$ then again  by the Mean Value Theorem (\ref{p22}) implies that there exists an $s$ between $\alpha$ and
$\sqrt{\alpha^2+ n^2\pi^2 \gamma^2\left(Q^{n,n}_{2,2}+2P^{n,n}_{1,2}\right)}$ such that
\bean 
P^{n,n}_{2,2}&=& {1\over 2  s^{1/2}}\left(Q^{n,n}_{2,2}+2P^{n,n}_{1,2}\right)\ \le \ {1\over 2  \alpha^{1/2}} {c\over n^r}
\eean
so again the sum (\ref{P22sum}) converges.

But because of $n^4 \pi^4 P^{n,n}_{2,2}$ term in (\ref{p22}) in order for  the sum 
\bea\label{P11sum}
P_{1,1}(x_1,x_2)&=&  \sum_{n=0}^\infty P^{n,n}_{1,1}\sin nx_1 \sin nx_2
\eea 
to converge $r$ must be larger than $8$ when $\alpha=0$ and $r$ must be larger than $5$ when $\alpha>0$. 

If  $\alpha>0$ then all of the closed loop eigenvalues (\ref{eigv2}) are in the open left half of the complex plane.
In particular 
for large $|n|$ the real parts of the closed loop eigenvalues (\ref{eigv2}) are more negative than ${\alpha\over2}$.

 Can we shift all the eigenvalues into the open left half of the complex plane if $\alpha=0$?  If $Q^{n,n}>0 $ but decays like ${1\over n^r}$
as $n\to \infty $ then the term outside the square root in (\ref{eigv2}) will be negative but it will decay in absolute value like ${1\over n^{r/2-1}}$.
For (\ref{P11sum}) to converge  $r$ must be greater than $8$ so the term outside the square root in (\ref{eigv2}) is converging to zero faster than ${1\over r^3}$.
But we are more interested in the convergence of feedback  (\ref{Ksum}) than the convergence of the optimal cost(\ref{Psum}).   
For  (\ref{Ksum}) to converge $r$ need only greater than $1$.   If we choose $1<r<2$ then the term outside the square root in (\ref{eigv2}) will grow
like $n^{2-r}$ so the higher the mode the higher the damping.  It is intersting to note that even if the optimal cost of an LQR problem does not exist, the LQR methodogy may yield a stabilizing feedbak.

\section{Conclusion}
We have used the simple and constructive technique of completing the square to solve the LQR problem for the stabilization of the linear beam equation using  boundary control.  The result  is an explicit formula for the quadratic optimal cost and the linear optimal feedback.  Our approach allows us to decouple the spatial frequencies so we can damp out all or just some frequencies.

\end{document}